\documentclass[10pt,a4paper]{article}
\usepackage{graphicx}
\usepackage[latin1]{inputenc}
\usepackage{amsmath}
\usepackage{amsfonts}
\usepackage{amssymb}
\usepackage{amsthm}
\usepackage{url}
\author{Alec Edgington \thanks{\tt alec@obtext.com}}
\title{Statistical regularities in the zeta zeros}
\date{December 2006}
\theoremstyle{plain} \newtheorem{conj}{Conjecture}
\theoremstyle{plain} 
\theoremstyle{plain} \newtheorem{thm}{Theorem}

\begin{document}

\maketitle

\begin{abstract}
Numerical investigations around a transformation of Landau's formula suggest certain statistical regularities in the distribution of zeros of the Riemann zeta function.
\end{abstract}

\section{Landau's formula}

Denote by $ Z $ the set of non-trivial zeros of the Riemann zeta function
\begin{equation*}
    Z = \{\rho \in \mathbb{C} : \zeta(\rho) = 0, 0 < \Re(\rho) < 1 \}
\end{equation*}
and let
\begin{equation*}
    S = \{-i(\rho-\frac{1}{2}) : \rho \in Z\}
\end{equation*}
(so that the Riemann hypothesis is just the assertion that $ S \subseteq \mathbb{R} $). The standard counting function for the zeros is
\begin{equation*}
    N(T) = \# ( S \cap [0,T] )
\end{equation*}
and its growth is
\begin{equation}\label{eq:N}
    N(T) = \frac{T}{2\pi} \log{\frac{T}{2\pi}} - \frac{T}{2\pi} + O(\log{T})
\end{equation}
as $ T \rightarrow \infty $.

Landau's formula is a statement about the asymptotic behaviour of the function
\begin{equation*}
    \lambda_a(T) = \sum_{\substack{t \in S \\ \lvert t \rvert \le T }} e^{iat}
\end{equation*}
for general $ a > 0 $. Notice that $ \lambda_a(T) $ is real, regardless of the Riemann hypothesis, since $ S = -S $ and any non-real $ t \in S $ come in conjugate pairs.

Landau\cite{land} proved the following theorem. Here $ \Lambda(x) $ is the von Mangoldt function, which is equal to $ \log{p} $ if $ x = p^k $ for $ p $ prime and $ k \in \mathbb{N}^+ $, and $ 0 $ otherwise.
\begin{thm}[Landau's formula] \label{thm:landau}
For all $ x > 1 $,
\begin{equation}\label{eq:lambda}
    \frac{\lambda_{\log{x}}(T)}{2T} = -\frac{1}{2\pi} \Lambda(x) x^{-\frac{1}{2}} + O(\frac{\log{T}}{T})
\end{equation}
as $ T \rightarrow \infty $.
\end{thm}
Thus $ \lambda_a(T) $ grows linearly with $ T $ when $ a = k \log{p} $. By contrast, when $ a \ne k \log{p} $ its growth is extremely slow: positive and negative values of $ \cos{(at)} $ for $ t \in S $ cancel out with great efficiency.

The ordinates of the zeros of the zeta function around height $ T $ share many statistical properties with the arguments of the eigenvalues of a random matrix from the circular unitary ensemble of order $ n $ where $ n \approx \log{\frac{T}{2\pi}} $. Now it is known\cite{diac} that the sum of the eigenvalues of such a matrix (all of which lie on the unit circle) is, in the limit as $ n \rightarrow \infty $, distributed as a Gaussian variable with \emph{fixed} variance. In other words, under the random-matrix model the sum of the cosines of a sequence of $ \log{\frac{T}{2\pi}} $ successive zeros around height $ T $ is statistically of order 1. But even this goes nowhere near explaining the actual degree of cancellation.

\section{Sum over a cycle}

Define
\begin{equation*}
    \eta_{a,h}(n) = \frac{1}{a} \log{\frac{n}{a}} + \sum_{\substack{t + h \in S \\ \lfloor \frac{at}{2\pi} \rfloor = n}} (e^{iat} - 1).
\end{equation*}
We have the following transformation of Landau's formula.
\begin{thm} \label{thm:eta}
If $ x > 1 $ and $ h \in \mathbb{R} $ then
\begin{equation*}
    \frac{1}{M} \sum_{0 \le n < M} \eta_{\log{x},h}(n) = -\frac{\Lambda(x)}{\log{x}} x^{-(\frac{1}{2}+ih)} + O(\frac{\log{M}}{M})
\end{equation*}
as $ M \rightarrow \infty $.
\end{thm}
\begin{proof}
Write $ a = \log{x} $ and $ T = \frac{2\pi}{a} M $. We use equations \ref{eq:N} and \ref{eq:lambda}, and Stirling's formula, to deduce
\begin{equation*}
\begin{split}
\sum_{0 \le n < M} \eta_{a,h}(n) & = \frac{1}{2} e^{-iah} \lambda_a(T) - N(T) + \frac{1}{a} \log{\Gamma{(\frac{aT}{2\pi})}} - \frac{T}{2\pi} \log{a} + O(\log{T}) \\
& = -\frac{T}{2\pi} \Lambda(x) x^{-(\frac{1}{2}+ih)} - \frac{T}{2\pi} \log{\frac{T}{2\pi}} + \frac{T}{2\pi} \\
& \qquad + \frac{T}{2\pi} \log{\frac{aT}{2\pi}} - \frac{T}{2\pi} - \frac{T}{2\pi} \log{a} + O(\log{T}) \\
& = -\frac{T}{2\pi} \Lambda(x) x^{-(\frac{1}{2}+ih)} + O(\log{T})
\end{split}
\end{equation*}
and the result follows.
\end{proof}

Numerical investigations, using data that Odlyzko has made available on his website\cite{odly}, suggest that the distribution of $ \eta_{a,h}(n) $ around its mean value is asymptotically independent of $ n $, as well as of $ h $.
\begin{conj}\label{conj:eta}
For every $ a > 0 $ there is a Borel function $ f_a : \mathbb{C} \rightarrow \mathbb{R} $ satisfying
\begin{equation*}
    \int_{\mathbb{C}} f_a = 1
\end{equation*}
such that for every $ h \in \mathbb{R} $ and every Borel set $ X \subseteq \mathbb{C} $,
\begin{equation*}
    \lim_{M \rightarrow \infty} \frac{\# \{ 0 \le n < M : \eta_{a,h}(n) + \frac{\Lambda(e^a)}{a} e^{-(\frac{1}{2} + ih)a} \in X \}}{M} = \int_X f_a.
\end{equation*}
\end{conj}

It is natural to consider the closely related continuous function $ H_a : \mathbb{R}^+ \rightarrow \mathbb{C} $ defined by
\begin{equation*}
    H_a(\tau) = \frac{1}{a} \log{\frac{\tau}{2\pi a}} - \sum_{\substack{u \in \tau-aS \\ \lvert u \rvert \le \pi}} (1 + e^{-iu}).
\end{equation*}
Here $ \tau-aS $ denotes $ \{ \tau-at : t \in S \} $. It is easily checked that for all $a$ and $h$,
\begin{equation} \label{eq:eta-H-relation}
    \eta_{a,h}(n) = H_a((2n+1)\pi + ah) + O(n^{-1})
\end{equation}
as $ n \rightarrow \infty $, so that the asymptotic distribution of $\eta_{a,h}(n)$ is determined by that of $H_a(\tau)$ for $ \tau \equiv \pi+ah $ modulo $2\pi$.

\begin{figure} [ht]
\centering
\includegraphics[width=120mm]{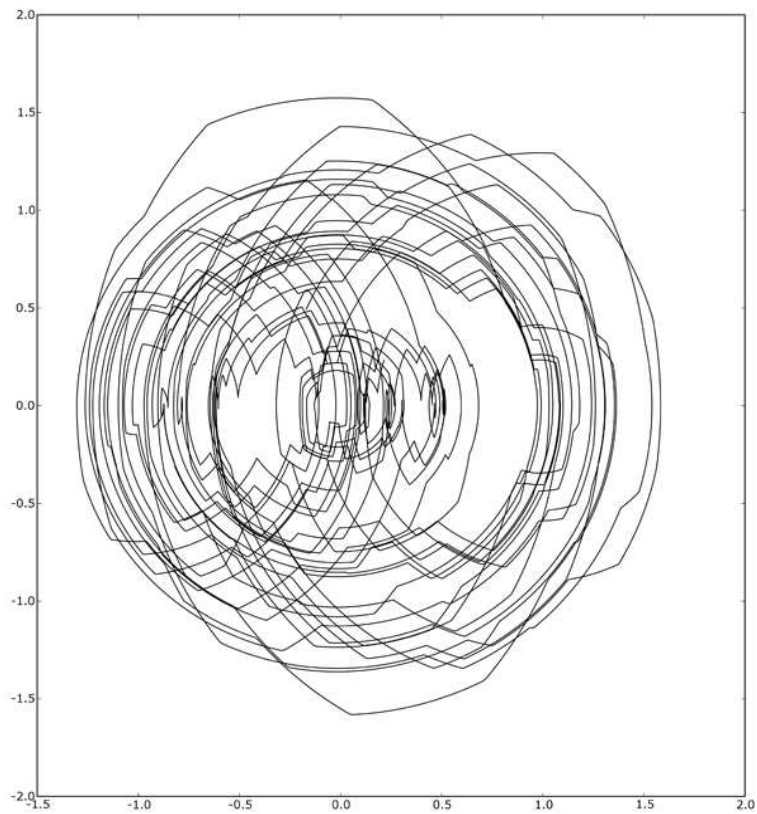}
\caption{\label{fig:H} Locus of $ H_1(\tau) $ for $ 50000 \le \tau \le 50600 $.}
\end{figure}

As $\tau$ increases continuously, $H_a(\tau)$ switches between near-circular orbits, as illustrated in Figure \ref{fig:H}. At `time' $\tau$ the orbit is centred on the point
\begin{displaymath}
\frac{1}{a} \log{\frac{\tau}{2\pi a}} - \nu_a(\tau)
\end{displaymath}
where
\begin{equation*}
    \nu_a(\tau) = \# ((\tau-aS) \cap [-\pi,\pi]).
\end{equation*}

These centres lie on the real axis at unit spacings, and they drift slowly in the positive direction. However, it seems that the statistical distribution of $H_a(\tau)$ is asymptotically independent of $\tau$.
\begin{conj}\label{conj:H}
For every $ a > 0 $ there is a Borel function $ F_a : \mathbb{C} \rightarrow \mathbb{R} $ satisfying
\begin{equation*}
    \int_{\mathbb{C}} F_a = 1
\end{equation*}
such that for every Borel set $ X \subseteq \mathbb{C} $,
\begin{equation*}
    \lim_{T \rightarrow \infty} \frac{\lambda ([0,T] \cap H_a^{-1}(X))}{T} = \int_X F_a
\end{equation*}
where $ \lambda $ denotes Lebesgue measure.
\end{conj}

Relation \ref{eq:eta-H-relation} implies that the function $F_a$ must be related to the function $f_a$ in Conjecture \ref{conj:eta} by
\begin{equation*}
    F_a(z) = \frac{1}{2\pi} \int_0^{2\pi} f_a(z - \frac{\Lambda(e^a)}{a} e^{-\frac{1}{2} a + i\theta}) \, \mathrm{d} \theta.
\end{equation*}
In particular, the two functions are the same when $ a \ne k \log{p} $.

\begin{figure} [ht]
\centering
\includegraphics[width=120mm]{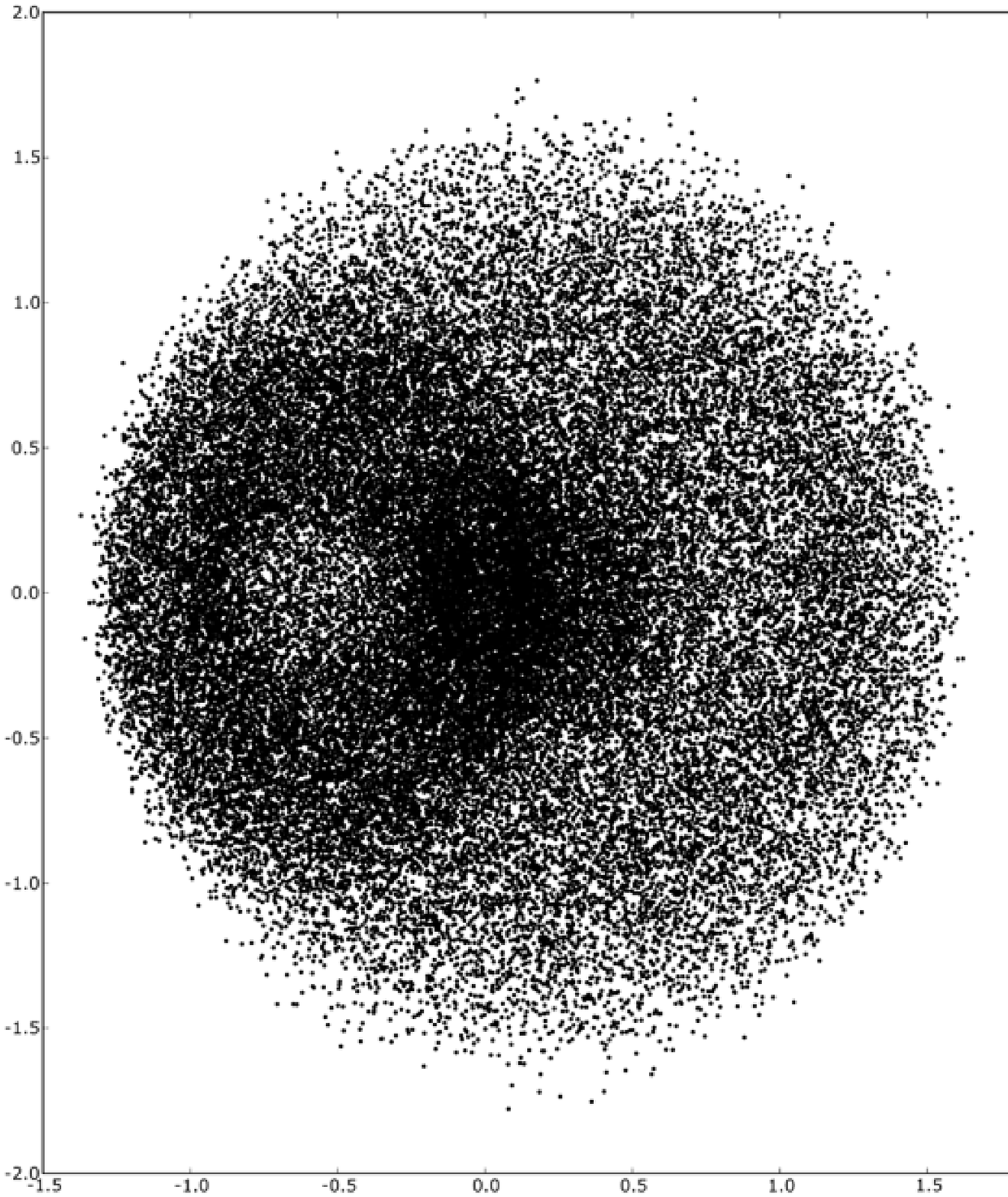}
\caption{\label{fig:H-dist-1} Values of $ H_1(\tau) $ in the complex plane, for $50000$ values of $\tau$ equally spaced between $\pi$ and $74920-\pi$.}
\end{figure}

\begin{figure} [ht]
\centering
\includegraphics[width=90mm]{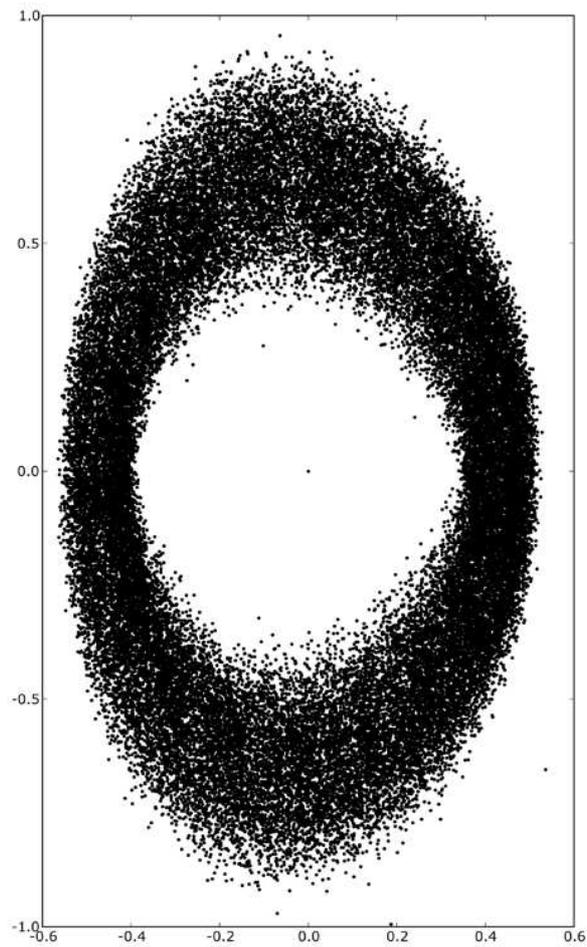}
\caption{\label{fig:H-dist-half} Values of $ H_{\frac{1}{2}}(\tau) $ in the complex plane, for $50000$ values of $\tau$ equally spaced between $\pi$ and $37460-\pi$.}
\end{figure}

\begin{figure} [ht]
\centering
\includegraphics[width=120mm]{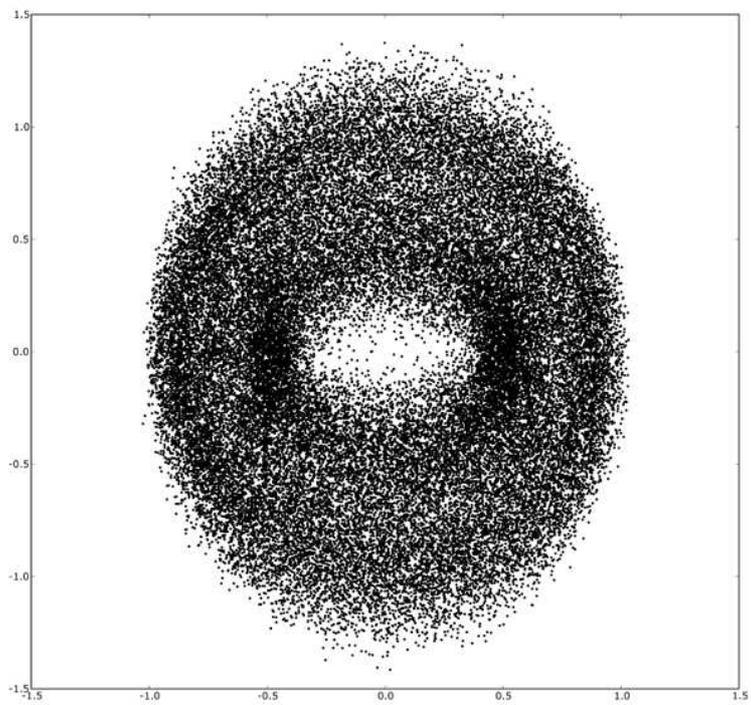}
\caption{\label{fig:H-dist-log2} Values of $ H_{\log{2}}(\tau) $ in the complex plane, for $50000$ values of $\tau$ equally spaced between $\pi$ and $74920\log{2} - \pi$.}
\end{figure}

\begin{figure} [ht]
\centering
\includegraphics[width=90mm]{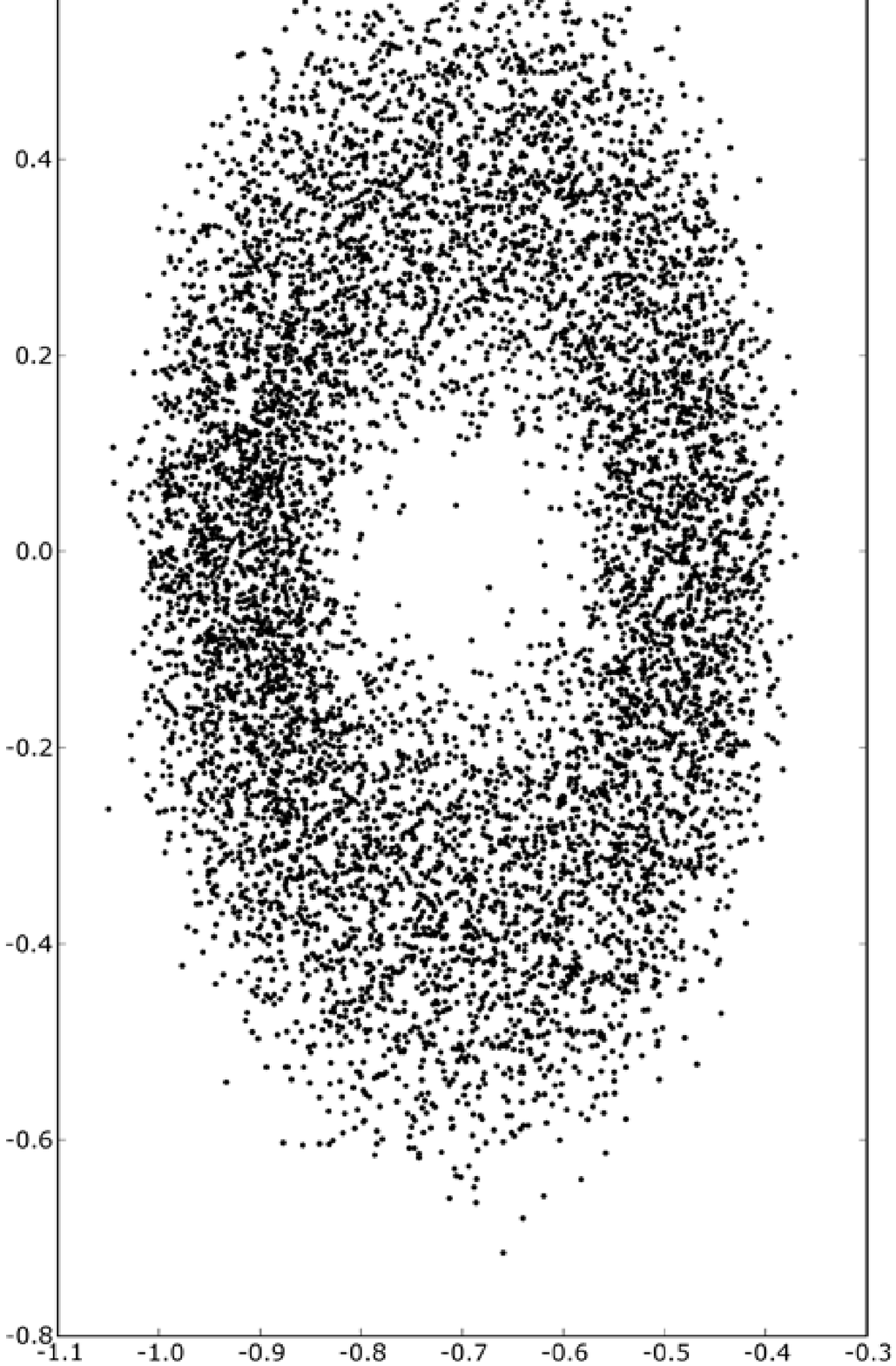}
\caption{\label{fig:eta-dist-log2} Values of $ \eta_{\log{2},0}(n) $ in the complex plane, for $ 4 \le n \le 8264 $.}
\end{figure}

Figures \ref{fig:H-dist-1}, \ref{fig:H-dist-half} and \ref{fig:H-dist-log2} illustrate the (empirical) density functions $F_a$ for $a=1$, $a=\frac{1}{2}$ and $a=\log{2}$ respectively. Figure \ref{fig:eta-dist-log2} illustrates $f_{\log{2}}$.

These plots are based on calculations from the first $100000$ zeros. The distribution seems to be independent of height, and plots based on zeros around the $ (10^{12})^{\textrm{th}} $, $ (10^{21})^{\textrm{st}} $ and $ (10^{22})^{\textrm{nd}} $ zero look similar. But even if the conjectures should in the end prove false, the patterns evident in these plots do demand an explanation.

\end{document}